\documentclass[a4paper]{article}
\usepackage{amsmath,amssymb,amsthm,latexsym}
\usepackage{graphicx,xcolor,tikz}
\renewenvironment{proof}{\noindent{\bf Proof.}}{~~$\Box$}
\theoremstyle{plain}
\newtheorem{thrm}{Theorem}[section]
\newtheorem{prot}[thrm]{Proposition}

\theoremstyle{definition}
\newtheorem{dfnt}[thrm]{Definition}%[section]
\newtheorem{remk}[thrm]{Remark}

\begin{document}
\title {On an upper bound of the set of copulas with a given curvilinear section} \vspace{8mm}
\author{ Yao Ouyang$^{a}$ \thanks{
         Email:oyy@zjhu.edu.cn} \qquad \qquad Yonghui Sun$^{a}$\qquad Hua-Peng Zhang$^b$\thanks{
         Corresponding author. Email:huapengzhang@163.com}\\
\small{\it $^a$ Faculty of Science, Huzhou University,
           Huzhou, Zhejiang 313000, China}\\
\small{\it $^b$ School of Science, Nanjing University of Posts and Telecommunications,
Nanjing 210023, China}}
\date{}
\maketitle
%
%%%%%%%%%%%%%%the text%%%%%%%%%%%%%%%%%%%%%%%%%%%%%%%%%%%%%%%%%%
\begin{abstract}
The characterizations when two natural upper bounds of the set of copulas with a given diagonal section are copulas have been well studied in the literature.
Given a curvilinear section, however, there is only a partial result concerning the characterization when a natural upper bound of the set of copulas is a copula.
In this paper, we completely solve the characterization problem for this natural upper bound to be a copula in the curvilinear case.

{\it Keywords} Copula; Diagonal section; Curvilinear section.
\end{abstract}
\section{Introduction}

(Bivariate) copulas, as functions joining bivariate distribution functions to their one-dimensional marginal distribution functions
according to Sklar's theorem~\cite{Skl59}, are indispensable for tackling the issues of non-normality and non-linearity in modern finance~\cite{CLV04}.
In view of the fact that the tail dependence of two random variables is determined by the diagonal and/or opposite diagonal sections of the underlying copula,
there have been several works devoted to studying copulas with given diagonal and/or opposite diagonal sections (see, for instance,~\cite{ADDF16, BMU09, DJ08, NQRU08}).
Given a diagonal section, while the smallest copula exists, the greatest one may not. This motivates the interest of characterizing when two natural upper bounds of
the set of copulas with a given diagonal section are copulas (see~\cite{KMS24, NQRU04, NQRU08, U08}).

By bending the diagonal into a curve, De Baets et al.~\cite{BMJ19} investigated the degree of asymmetry of a (quasi-)copula with respect to a curve. In particular,
they obtained the smallest copula with a given curvilinear section. Zou et al.~\cite{ZSX22} continued to present the greatest quasi-copula with a given curvilinear section.
Under some additional assumptions, they characterized when the greatest quasi-copula with a given curvilinear section is a copula. The aim of this paper is to completely solve the characterization when the greatest quasi-copula with a given curvilinear section is a copula without any additional assumption.

The remainder of this paper is organized as follows. We recall some basic notions and results related to (quasi-)copulas as well as that with a given curvilinear section
in Section~\ref{Sec-Preliminaries}. Section~\ref{Sec-Main} is devoted to the characterization theorem. We end with future work in Section~\ref{Sec-Conclu}.

\section{Preliminaries}\label{Sec-Preliminaries}

In this section, we recall some basic notions and results related to (quasi-)copulas as well as that with a given curvilinear section.
\begin{dfnt}\cite{Nel06}
A (bivariate) \emph{copula} is a function $C: [0, 1]^2\to [0, 1]$ satisfying:\\
(C1) the boundary condition, \emph{i.e.},
 \[
 C(u, 0)=C(0, u)=0\ {\rm and} \ C(u, 1)=C(1, u)=u\ (\forall\,u\in [0, 1]);
 \]
(C2) the $2$-increasing property, \emph{i.e.},
$V(R)\geq 0$ for any rectangle $R=[u_1, u_2]\times [v_1, v_2]$,
where $V(R)=C(u_2, v_2)-C(u_2, v_1)-C(u_1, v_2)+C(u_1, v_1)$.
\end{dfnt}

\begin{dfnt}\cite{Nel06}
A (bivariate) \emph{quasi-copula} is a function $Q: [0, 1]^2\to [0, 1]$ satisfying (C1) and the following conditions:\\
(Q1) $Q$ is increasing in each variable;\\
(Q2) $Q$ is Lipschitz, \emph{i.e.}, for all $u_1, u_2, v_1, v_2\in [0, 1]$, it holds that
\[
|Q(u_1, v_1)-Q(u_2, v_2)|\leq |u_1-u_2|+|v_1-v_2|.
\]
\end{dfnt}
Note that (C1) and (C2) together imply both (Q1) and (Q2), \emph{i.e.}, every copula is a quasi-copula, but not vice versa.
Three important copulas are $W, M$ and $\Pi$ given by
\[W(u, v)=\max\{u+v-1, 0\}, M(u, v)=\min\{u, v\} \ {\rm and} \ \Pi(u, v)=uv.\]
$W$ and $M$ are known as the \emph{Fr\'{e}chet-Hoeffding bounds} for copulas and quasi-copulas since $W\leq Q\leq M$
for any quasi-copula $Q$, while $\Pi$ is known as the independence copula.

\begin{dfnt}\cite{NQRU04}
A \emph{diagonal function} is a function $\delta: [0, 1]\to [0, 1]$ with the properties:\\
(D1) $\delta(1)=1$ and $\delta(t)\leq t$ for all $t\in [0, 1]$;\\
(D2) $0\leq \delta(t_2)-\delta(t_1)\leq 2(t_2-t_1)$ for all $t_1, t_2\in [0, 1]$ with $t_1\leq t_2$.
\end{dfnt}

The diagonal section $\delta_Q: [0, 1]\to [0, 1]$ of any quasi-copula $Q$ defined by $\delta_Q(t)=Q(t, t)$ is a diagonal function. Conversely, for any diagonal function $\delta$, there exist quasi-copulas with diagonal section $\delta$,
and the smallest quasi-copula $B_\delta$ with diagonal section $\delta$ is given by
\[
B_\delta(u, v) =\left \{
        \begin {array}{ll}
         v-\min\limits_{t\in[v, u]}\{t-\delta(t)\}
                  &\quad \text{if\ \ $v\leq u$}  \\[2mm]
           u-\min\limits_{t\in[u, v]}\{t-\delta(t)\}
                  &\quad \text{if\ \ $v>u$,}
        \end {array}
       \right.
\]
while the greatest quasi-copula $A_\delta$ with diagonal section $\delta$ is given by
\[
A_\delta(u, v) =\left \{
        \begin {array}{ll}
         \min\left\{v, u-\max\limits_{t\in[v, u]}\{t-\delta(t)\}\right\}
                  &\quad \text{if\ \ $v\leq u$}  \\[2mm]
           \min\left\{u, v-\max\limits_{t\in[u, v]}\{t-\delta(t)\}\right\}
                  &\quad \text{if\ \ $v>u$.}
        \end {array}
       \right.
\]
Interestingly, $B_\delta$ is also a copula (called the \emph{Bertino copula}), but $A_\delta$ may not be.
There have been characterization theorems for $A_\delta$ to be a copula (see, for instance, \cite{NQRU04, U08}).

Let $\Phi$ be the set of all functions $\phi: [0, 1]\to [0, 1]$ that is continuous and strictly increasing with $\phi(0)=0$ and $\phi(1)=1$.
Obviously, ${\rm id}\in\Phi$, where ${\rm id}$ is the identity map on $[0, 1]$, \emph{i.e.}, ${\rm id}(x)=x$ for all $x\in [0, 1]$.
For any $\phi\in\Phi$, the curvilinear section $g_{Q, \phi}: [0, 1]\to [0, 1]$ of a quasi-copula $Q$ defined by $g_{Q, \phi}(t)=Q(t, \phi(t))$
has properties:\\
(i) $\max\{0, t+\phi(t)-1\}\leq g_{Q,\phi}(t)\leq\min\{t, \phi(t)\}$ for all $t\in [0, 1]$;\\
(ii) $0\leq g_{Q,\phi}(t_2)-g_{Q,\phi}(t_1)\leq t_2-t_1+\phi(t_2)-\phi(t_1)$ for all $t_1, t_2\in [0, 1]$ with $t_1<t_2$.\\
Note that property (ii) can be rewritten as: $g_{Q, \phi}(t)$ is increasing and $g_{Q, \phi}(t)-t-\phi(t)$ is decreasing.
For any $\phi\in\Phi$, denote by $\Delta_\phi$ the set of all functions $\Gamma_\phi: [0, 1]\to [0, 1]$ with properties (i) and (ii).
Particularly, if $\phi={\rm id}$, then properties (i) and (ii) are equivalent to (D1) and (D2), and thus $\Delta_{\rm id}$ is nothing else but the set of all diagonal functions.

For any $\phi\in\Phi$ and $\Gamma_\phi\in\Delta_\phi$, similarly to the diagonal case, there exist quasi-copulas with curvilinear section $\Gamma_\phi$,
and the smallest quasi-copula $B_{\Gamma_\phi}$ with curvilinear section $\Gamma_\phi$ is given by

\begin{equation}\label{Eq-DfntofBGamma}
B_{\Gamma_\phi}(u, v) =\left \{
        \begin {array}{ll}
         v-\min\limits_{t\in[\phi^{-1}(v), u]}\{\phi(t)-\Gamma_\phi(t)\}
                  &\quad \text{if\ \ $v\leq\phi(u)$}  \\[2mm]
           u-\min\limits_{t\in[u, \phi^{-1}(v)]}\{t-\Gamma_\phi(t)\}
                  &\quad \text{if\ \ $v>\phi(u)$,}
        \end {array}
       \right.
\end{equation}
while the greatest quasi-copula $A_{\Gamma_\phi}$ with curvilinear section $\Gamma_\phi$ is given by
\begin{equation}\label{Eq-DfntofAGamma}
A_{\Gamma_\phi}(u, v) =\left \{
        \begin {array}{ll}
         \min\left\{v, u-\max\limits_{t\in[\phi^{-1}(v), u]}\{t-\Gamma_\phi(t)\}\right\}
                  & \text{if\ \ $v\leq\phi(u)$}  \\[2mm]
           \min\left\{u, v-\max\limits_{t\in[u, \phi^{-1}(v)]}\{\phi(t)-\Gamma_\phi(t)\}\right\}
                  & \text{if\ \ $v>\phi(u)$.}
        \end {array}
       \right.
\end{equation}
$B_{\Gamma_\phi}$ is also a copula, called the \emph{Bertino copula with curvilinear section} $\Gamma_\phi$~\cite{BMJ19}.
However, $A_{\Gamma_\phi}$ may not be a copula. Zou et al.~\cite{ZSX22} characterized when $A_{\Gamma_\phi}$ is a copula under some additional assumptions.

Note that $m_\phi\in\Delta_\phi$ for all $\phi\in\Phi$, where $m_\phi: [0, 1]\to [0, 1]$ is defined by $m_\phi(t)=\min\{t, \phi(t)\}$.
It is easy to check that $A_{m_\phi}=M$. If $\phi={\rm id}$, then $B_{m_\phi}=M$. If $\phi\not={\rm id}$, then $B_{m_\phi}\not=M$, as
the following proposition shows.
\begin{prot}
Let $\phi\in\Phi$. Then $B_{m_\phi}=M$ if and only if $\phi={\rm id}$.
\end{prot}
\begin{proof}
We only need to prove the necessity. Suppose that $\phi\not=id$, \emph{i.e.}, there exists $u\in \,]0, 1[\,$ such that $\phi(u)\not= u$.
We consider the case $\phi(u)<u$ (the proof is similar for the case $\phi(u)>u$). It follows from the continuity of $\phi$ that there exists $v\in\,]u, 1[\,$
such that $\phi(v)<u$. Hence,
\[
B_{m_\phi}(u, \phi(v))=u-\min_{t\in [u, v]}\{t-\phi(t)\}.
\]
Since $u-(t-\phi(t))=(u-t)+\phi(t)<\phi(v)$ for all $t\in [u, v]$, we have
\[
B_{m_\phi}(u, \phi(v))<\phi(v)=M(u, \phi(v)),
\]
a contradiction.
\end{proof}

Let $\phi\in\Phi$ and $\{\,]a_i, b_i[\,\}_{i\in I}$ be pairwisely disjoint countably many open subintervals of $\,]0, 1[\,$.
We say that $\Gamma_\phi: [0, 1]\to [0, 1]$ is determined by $\{\,]a_i, b_i[\,\}_{i\in I}$ provided that
\[
\Gamma_\phi(t) =\left \{
        \begin {array}{ll}
         m_\phi(a_i)
                  &\quad \text{\rm if\ \ $t\in \,]a_i, u_i^\ast]$} \\[2mm]
           \phi(t)+t-\max\{b_i, \phi(b_i)\}
                  &\quad \text{\rm if\ \ $t\in [u_i^\ast, b_i[\,$} \\[2mm]
         m_\phi(t)
                  &\quad {\rm otherwise,}
        \end {array}
       \right.
\]
where $u_i^\ast\in\,]a_i, b_i[\,$ is the unique solution to the equation:
$\phi(t)+t=m_\phi(a_i)+\max\{b_i, \phi(b_i)\}$, \emph{i.e.}, $\phi(u_i^\ast)+u_i^\ast=m_\phi(a_i)+\max\{b_i, \phi(b_i)\}$.
It is easy to prove that $\Gamma_\phi\in\Delta_\phi$. As we shall see in next section, $A_{\Gamma_\phi}$ is a copula only if
$\Gamma_\phi=m_\phi$ or $\Gamma_\phi$ is determined by pairwisely disjoint countably many open subintervals of $\,]0, 1[\,$.

\section{Characterization theorem}\label{Sec-Main}

In this section, we characterize when $A_{\Gamma_\phi}$ is a copula. To simplify the proof of the necessity of the characterization theorem,
we need the help of the following two propositions.
\begin{prot}\label{prot-featureofphi}
Let $A_{\Gamma_\phi}$ be a copula and $a\in\,]0, 1[\,$. If $\Gamma_\phi(a)=m_\phi(a)$, then $b_0\geq\max\{\phi(a_0), \phi^{-1}(a_0)\}$, where
$b_0=\inf\{t\in\,]a, 1[\,\mid \Gamma_\phi(t)<m_\phi(t)\}$ and $a_0=\sup\{t\in\,]0, a[\,\mid \Gamma_\phi(t)<m_\phi(t)\}$, stipulating $\inf\emptyset=1$ and $\sup\emptyset=0$.
\end{prot}
\begin{proof}
Suppose that $b_0<\phi(a_0)$ or $b_0<\phi^{-1}(a_0)$. We only consider the case $b_0<\phi(a_0)$ (the proof is similar for the other case). Obviously, $0<a_0\leq a\leq b_0<1$.
It follows from the continuity of $\phi(t)$ at $a_0$ and the definitions of $a_0, b_0$ that there exist $a_1\in\,]0, a_0[\,$ and $b_1\in\,]b_0, 1[\,$ such that $b_1<\phi(a_1)$, $\Gamma_\phi(a_1)<m_\phi(a_1)$ and $\Gamma_\phi(b_1)<m_\phi(b_1)$.
For the rectangle $R=[a, b_1]\times [\phi(a_1), \phi(a)]$, we have
\[A_{\Gamma_\phi}(a, \phi(a_1))=\min\{\phi(a_1), a-\max\limits_{t\in [a_1, a]}\{t-\Gamma_\phi(t)\}\}
=a-\max\limits_{t\in [a_1, a]}\{t-\Gamma_\phi(t)\},\]
\[A_{\Gamma_\phi}(a, \phi(a))=\Gamma_\phi(a)=m_\phi(a)=a,\]
\[A_{\Gamma_\phi}(b_1, \phi(a_1))=\min\{\phi(a_1), b_1-\max\limits_{t\in [a_1, b_1]}\{t-\Gamma_\phi(t)\}\}
=b_1-\max\limits_{t\in [a_1, b_1]}\{t-\Gamma_\phi(t)\},\]
and
\[A_{\Gamma_\phi}(b_1, \phi(a))=\min\{\phi(a), b_1-\max\limits_{t\in [a, b_1]}\{t-\Gamma_\phi(t)\}\}
=b_1-\max\limits_{t\in [a, b_1]}\{t-\Gamma_\phi(t)\}.\]
Hence,
\begin{eqnarray*}
V(R)= \max\limits_{t\in [a_1, b_1]}\{t-\Gamma_\phi(t)\}-\max\limits_{t\in [a_1, a]}\{t-\Gamma_\phi(t)\}-\max\limits_{t\in [a, b_1]}\{t-\Gamma_\phi(t)\}.
\end{eqnarray*}
Note that $\max\limits_{t\in [a_1, a]}\{t-\Gamma_\phi(t)\}\geq a_1-\Gamma_\phi(a_1)>a_1-m_\phi(a_1)=0$. Similarly, $\max\limits_{t\in [a, b_1]}\{t-\Gamma_\phi(t)\}>0$.
Therefore, we have $V(R)<0$, a contradiction.
\end{proof}

\begin{remk}\label{remk-featureofphi}
For all $a, b\in[0, 1]$ with $a\leq b$, it is easy to see that $b\geq\max\{\phi(a), \phi^{-1}(a)\}$ if and only if $m_\phi(b)\geq\max\{a, \phi(a)\}$.
Hence, under the assumptions in Proposition~\ref{prot-featureofphi}, we have $\phi(a)=a$ when $a_0=b_0$.
\end{remk}

\begin{prot}\label{prot-featureofGamma}
Let $A_{\Gamma_\phi}$ be a copula and $\Gamma_\phi(t)<m_\phi(t)$ for all $t\in\,]a, b[\,$, where $a, b\in[0, 1]$ with $a<b$. Then the following results hold:
\begin{itemize}%\setlength{\itemindent}{5pt}
\item[{\rm(i)}] $\Gamma_\phi(t)$ is strictly increasing on $[a^\ast, b]$, where
\[a^\ast=\max\{t\in [a, b]\mid \Gamma_\phi(t)=\Gamma_\phi(a)\};\]
\item[{\rm(ii)}] $\Gamma_\phi(t)-t-\phi(t)$ is strictly decreasing on $[a, b^\ast]$, where
\[b^\ast=\min\{t\in [a, b]\mid \Gamma_\phi(t)-t-\phi(t)=\Gamma_\phi(b)-b-\phi(b)\};\]
\item[{\rm(iii)}] $a^\ast=b^\ast$;
\item[{\rm(iv)}] if $\Gamma_\phi(a)-m_\phi(a)=\Gamma_\phi(b)-m_\phi(b)=0$, then $a^\ast\in\,]a, b[\,$ and is the unique solution to the equation:
$\phi(t)+t=m_\phi(a)+\max\{b, \phi(b)\}$, \emph{i.e.}, $\phi(a^\ast)+a^\ast=m_\phi(a)+\max\{b, \phi(b)\}$, and
\[
\Gamma_\phi(t) =\left \{
        \begin {array}{ll}
         m_\phi(a)
                  &\quad \text{\rm{if}\ \ $t\in [a, a^\ast]$}  \\[2mm]
           \phi(t)+t-\max\{b, \phi(b)\}
                  &\quad \text{\rm{if}\ \ $t\in[a^\ast, b]$.}
        \end {array}
       \right.
\]
\end{itemize}
\end{prot}
\begin{proof}
(i) Suppose that there exist $b_1, b_2\in\,]a^\ast, b]$ with $b_1<b_2$ such that $\Gamma_\phi (b_1)=\Gamma_\phi(b_2)$.
Set $b_0=\min\{t\in [a^\ast, b_2[\,\mid \Gamma_\phi(t)=\Gamma_\phi(b_2)\}$. Obviously, $a^\ast<b_0\leq b_1$ and
$\Gamma_\phi(t)<\Gamma_\phi(b_0)$ for all $t<b_0$. Since $\Gamma_\phi(b_0)<m_\phi(b_0)$ and $m_\phi(t)$ is continuous, there exists
$a_0\in\,]a^\ast, b_0[\,$ such that $\Gamma_\phi(b_0)<m_\phi(a_0)$. For the rectangle $R=[a_0, b_0]\times [\phi(a_0), \phi(b_2)]$, we have
\[A_{\Gamma_\phi}(a_0, \phi(b_2))=\min\{a_0, \phi(b_2)-\max\limits_{t\in [a_0, b_2]}\{\phi(t)-\Gamma_\phi(t)\}\}
=\phi(b_2)-\max\limits_{t\in [a_0, b_2]}\{\phi(t)-\Gamma_\phi(t)\},\]
\[A_{\Gamma_\phi}(b_0, \phi(b_2))=\min\{b_0, \phi(b_2)-\max\limits_{t\in [b_0, b_2]}\{\phi(t)-\Gamma_\phi(t)\}\}
=\phi(b_2)-\max\limits_{t\in [b_0, b_2]}\{\phi(t)-\Gamma_\phi(t)\},\]
and
\[A_{\Gamma_\phi}(b_0, \phi(a_0))=\min\{\phi(a_0), b_0-\max\limits_{t\in [a_0, b_0]}\{t-\Gamma_\phi(t)\}\}
=b_0-\max\limits_{t\in [a_0, b_0]}\{t-\Gamma_\phi(t)\}.\]
Hence,
\begin{eqnarray*}
V(R)= \Gamma_\phi(a_0)-b_0+\max\limits_{t\in [a_0, b_0]}\{t-\Gamma_\phi(t)\}+\max\limits_{t\in [a_0, b_2]}\{\phi(t)-
\Gamma_\phi(t)\}-\max\limits_{t\in [b_0, b_2]}\{\phi(t)-\Gamma_\phi(t)\}.
\end{eqnarray*}
We distinguish two cases to get contradictions.
\begin{itemize}%\setlength{\itemindent}{5pt}
\item[{\rm(a)}] $\max\limits_{t\in [a_0, b_2]}\{\phi(t)-\Gamma_\phi(t)\}=\max\limits_{t\in [b_0, b_2]}\{\phi(t)-\Gamma_\phi(t)\}$.

In this case, we have $V(R)= \Gamma_\phi(a_0)-b_0+\max\limits_{t\in [a_0, b_0]}\{t-\Gamma_\phi(t)\}$. Since $\Gamma_\phi(a_0)<\Gamma_\phi(b_0)$,
it holds that \[\Gamma_\phi(a_0)-b_0+t-\Gamma_\phi(t)=(\Gamma_\phi(a_0)-\Gamma_\phi(t))+(t-b_0)<0\]
for all $t\in [a_0, b_0]$, which implies that $V(R)<0$,
a contradiction.

\item[{\rm(b)}] $\max\limits_{t\in [a_0, b_2]}\{\phi(t)-\Gamma_\phi(t)\}=\max\limits_{t\in [a_0, b_0]}\{\phi(t)-\Gamma_\phi(t)\}$.

In this case, we have
\begin{eqnarray*}
V(R)= \Gamma_\phi(a_0)-b_0+\max\limits_{t\in [a_0, b_0]}\{t-\Gamma_\phi(t)\}+\max\limits_{t\in [a_0, b_0]}\{\phi(t)-
\Gamma_\phi(t)\}-\max\limits_{t\in [b_0, b_2]}\{\phi(t)-\Gamma_\phi(t)\}.
\end{eqnarray*}
Since $\Gamma_\phi(t)$ takes constant value $\Gamma_\phi(b_0)$ on $[b_0, b_2]$, $\max\limits_{t\in [b_0, b_2]}\{\phi(t)-\Gamma_\phi(t)\}=\phi(b_2)-\Gamma_\phi(b_0)$.
Hence, \begin{eqnarray*}
V(R)= \Gamma_\phi(a_0)-b_0+\max\limits_{t\in [a_0, b_0]}\{t-\Gamma_\phi(t)\}+\max\limits_{t\in [a_0, b_0]}\{\phi(t)-
\Gamma_\phi(t)\}-\phi(b_2)+\Gamma_\phi(b_0).
\end{eqnarray*}
For all $s, t\in [a_0, b_0]$, setting $m=\min\{s, t\}$ and $M=\max\{s, t\}$, it holds that
\begin{eqnarray*}
&  & \Gamma_\phi(a_0)-b_0+s-\Gamma_\phi(s)+\phi(t)-\Gamma_\phi(t)-\phi(b_2)+\Gamma_\phi(b_0)\\
&= & \Gamma_\phi(a_0)-b_0+(s-M)+M-\Gamma_\phi(m)+(\phi(t)-\phi(M))+\phi(M)\\
& & -\Gamma_\phi(M)-\phi(b_0)+(\phi(b_0)-\phi(b_2))+\Gamma_\phi(b_0)\\
&= &  (\Gamma_\phi(a_0)-\Gamma_\phi(m))+((\Gamma_\phi(b_0)-b_0-\phi(b_0))-(\Gamma_\phi(M)-M-\phi(M)))\\
& & +(s-M)+(\phi(t)-\phi(M))+(\phi(b_0)-\phi(b_2))\\
&< &  0.
\end{eqnarray*}
Therefore, we have $V(R)<0$, a contradiction.
\end{itemize}

(ii) Suppose that there exist $a_1, a_2\in [a, b^\ast[\,$ with $a_1<a_2$ such that $\Gamma_\phi(a_1)-a_1-\phi(a_1)=\Gamma_\phi(a_2)-a_2-\phi(a_2)$,
which implies that $\Gamma_\phi(a_1)<\Gamma_\phi(a_2)$, and so $a_2>a^\ast$.
Set $a_0=\max\{t\in \,]a_1, b^\ast]\mid \Gamma_\phi(t)-t-\phi(t)=\Gamma_\phi(a_1)-a_1-\phi(a_1)\}$. Obviously, $a^\ast<a_2\leq a_0<b^\ast$,
and $\Gamma_\phi(t)>\Gamma_\phi(a_0)$ (by (i)) and $\Gamma_\phi(t)-t-\phi(t)<\Gamma_\phi(a_0)-a_0-\phi(a_0)$ for all $t>a_0$.
Since $\Gamma_\phi(a_0)<m_\phi(a_0)$ and $\Gamma_\phi(t)$ is continuous, there exists
$b_0\in\,]a_0, b^\ast[\,$ such that $\Gamma_\phi(b_0)<m_\phi(a_0)$. For the rectangle $R=[a_0, b_0]\times [\phi(a_0), \phi(b_0)]$, we have
\[A_{\Gamma_\phi}(a_0, \phi(b_0))=\min\{a_0, \phi(b_0)-\max\limits_{t\in [a_0, b_0]}\{\phi(t)-\Gamma_\phi(t)\}\}
=\phi(b_0)-\max\limits_{t\in [a_0, b_0]}\{\phi(t)-\Gamma_\phi(t)\}\]
and
\[A_{\Gamma_\phi}(b_0, \phi(a_0))=\min\{\phi(a_0), b_0-\max\limits_{t\in [a_0, b_0]}\{t-\Gamma_\phi(t)\}\}
=b_0-\max\limits_{t\in [a_0, b_0]}\{t-\Gamma_\phi(t)\}.\]
Hence,
\begin{eqnarray*}
V(R)= \Gamma_\phi(b_0)-b_0-\phi(b_0)+\Gamma_\phi(a_0)+\max\limits_{t\in [a_0, b_0]}\{t-\Gamma_\phi(t)\}+\max\limits_{t\in [a_0, b_0]}\{\phi(t)-
\Gamma_\phi(t)\}.
\end{eqnarray*}
For all $s, t\in [a_0, b_0]$, setting $m=\min\{s, t\}$ and $M=\max\{s, t\}$, it holds that
\begin{eqnarray*}
&  & \Gamma_\phi(b_0)-b_0-\phi(b_0)+\Gamma_\phi(a_0)+s-\Gamma_\phi(s)+\phi(t)-\Gamma_\phi(t)\\
&= & \Gamma_\phi(b_0)-b_0-\phi(b_0)+\Gamma_\phi(a_0)+(s-M)+M-\Gamma_\phi(m)+(\phi(t)-\phi(M))+\phi(M)-\Gamma_\phi(M)\\
&= & ((\Gamma_\phi(b_0)-b_0-\phi(b_0))-(\Gamma_\phi(M)-M-\phi(M)))+(\Gamma_\phi(a_0)-\Gamma_\phi(m))\\
&  &  +(s-M)+(\phi(t)-\phi(M))\\
&\leq  & 0,
\end{eqnarray*}
where the last equality holds only if $s=t=M$, and so $m=M$. It follows from $\Gamma_\phi(t)>\Gamma_\phi(a_0)$ for all $t>a_0$ and $\Gamma_\phi(b_0)-b_0-\phi(b_0)<\Gamma_\phi(a_0)-a_0-\phi(a_0)$ that the last equality cannot hold. Therefore, we have $V(R)<0$, a contradiction.

(iii) It is easy to see that $\Gamma_\phi(t)$ is constant on $[a, a^\ast]$ and is strictly increasing on $[b^\ast, b]$, and thus $a^\ast\leq b^\ast$.
Suppose that $a^\ast< b^\ast$. Then there exist $a_0, b_0\in\,]a^\ast, b^\ast[\,$ with $a_0< b_0$ such that $\Gamma_\phi(b_0)<m_\phi(a_0)$. In a similar way as (ii),
for the rectangle $R=[a_0, b_0]\times [\phi(a_0), \phi(b_0)]$, we have $V(R)<0$, a contradiction.

(iv) It follows from (i)-(iii) that
\[
\Gamma_\phi(t) =\left \{
        \begin {array}{ll}
         m_\phi(a)
                  &\quad \text{\rm{if}\ \ $t\in [a, a^\ast]$}  \\[2mm]
           \phi(t)+t-\max\{b, \phi(b)\}
                  &\quad \text{\rm{if}\ \ $t\in[a^\ast, b]$.}
        \end {array}
       \right.
\]
The well-definedness of $\Gamma_\phi(a^\ast)$ implies that $\phi(a^\ast)+a^\ast=m_\phi(a)+\max\{b, \phi(b)\}$. Clearly, $a^\ast\in\,]a, b[\,$ and the solution to the equation $\phi(t)+t=m_\phi(a)+\max\{b, \phi(b)\}$ is unique.
\end{proof}

Now, we are ready to present the characterization theorem.
\begin{thrm}\label{Thrm-main}
Let $\phi\in\Phi$, $\Gamma_\phi\in\Delta_\phi$ and $A_{\Gamma_\phi}$ be defined as in~\eqref{Eq-DfntofAGamma}. Then $A_{\Gamma_\phi}$ is a copula if and only if
$\Gamma_\phi=m_\phi$ or there exist pairwisely disjoint countably many open subintervals $\{\,]a_i, b_i[\,\}_{i\in I}$ of $\,]0, 1[\,$ such that
\begin{itemize}%\setlength{\itemindent}{5pt}
\item[\rm (i)] $\Gamma_\phi$ is determined by $\{\,]a_i, b_i[\,\}_{i\in I}$, \emph{i.e.},
\[
\Gamma_\phi(t) =\left \{
        \begin {array}{ll}
         m_\phi(a_i)
                  &\quad \text{\rm if\ \ $t\in \,]a_i, u_i^\ast]$} \\[2mm]
           \phi(t)+t-\max\{b_i, \phi(b_i)\}
                  &\quad \text{\rm if\ \ $t\in [u_i^\ast, b_i[\,$} \\[2mm]
         m_\phi(t)
                  &\quad {\rm otherwise,}
        \end {array}
       \right.
\]
where $u_i^\ast\in\,]a_i, b_i[\,$ is the unique solution to the equation:
$\phi(t)+t=m_\phi(a_i)+\max\{b_i, \phi(b_i)\}$, \emph{i.e.}, $\phi(u_i^\ast)+u_i^\ast=m_\phi(a_i)+\max\{b_i, \phi(b_i)\}$;

\item[\rm (ii)] $\phi$ is compatible with $\{\,]a_i, b_i[\,\}_{i\in I}$, \emph{i.e.}, $m_\phi(a_j)\geq\max\{b_i, \phi(b_i)\}$ for all $i, j\in I$ with $a_j\geq b_i$.
\end{itemize}
\end{thrm}
\begin{proof}
\emph{Necessity:} Suppose that $\Gamma_\phi\not=m_\phi$. Since $\Gamma_\phi(0)=m_\phi(0)$, $\Gamma_\phi(1)=m_\phi(1)$, and $\Gamma_\phi$ and $m_\phi$ are continuous,
$\{t\in [0, 1] \mid m_\phi(t)-\Gamma_\phi(t)>0\}$ is a nonempty open set. Hence, there exist pairwisely disjoint countably many open subintervals $\{\,]a_i, b_i[\,\}_{i\in I}$ of $\,]0, 1[\,$ such that $\{t\in [0, 1] \mid m_\phi(t)-\Gamma_\phi(t)>0\}=\bigcup\limits_{i\in I}\,]a_i, b_i[\,$. It follows from Proposition~\ref{prot-featureofGamma}(iv) that (i) holds, while (ii) follows from Proposition~\ref{prot-featureofphi} and Remark~\ref{remk-featureofphi}.

\emph{Sufficiency:} If $\Gamma_\phi=m_\phi$, then $A_{\Gamma_\phi}=M$ is a copula.
Suppose that there exist pairwisely disjoint countably many open subintervals $\{\,]a_i, b_i[\,\}_{i\in I}$ of $\,]0, 1[\,$ such that
$\Gamma_\phi$ is determined by $\{\,]a_i, b_i[\,\}_{i\in I}$ and $\phi$ is compatible with $\{\,]a_i, b_i[\,\}_{i\in I}$.
Since $A_{\Gamma_\phi}$ is a quasi-copula (see Proposition 12 of~\cite{ZSX22}), it remains to show that $A_{\Gamma_\phi}$ is $2$-increasing.
Note that any rectangle $R=[a, b]\times [c, d]$ can be partitioned into at most three non-overlapping rectangles $R_i$ and $V(R)$ is the sum of all $V(R_i)$s,
where each $R_i$ belongs to one of the three types of rectangles:

\noindent Type 1: $R$ is above the curve of $\phi$, \emph{i.e.}, $R\subset\{(u, v)\mid v\geq\phi(u)\}$;

\noindent Type 2: $R$ is below the curve of $\phi$, \emph{i.e.}, $R\subset\{(u, v)\mid v\leq\phi(u)\}$;

\noindent Type 3: $R$ has two vertices on the curve of $\phi$, \emph{i.e.}, $R=[u_1, u_2]\times [\phi(u_1), \phi(u_2)]$.

\noindent We only need to prove that $V(R)\geq 0$ for each of the three types of rectangles $R=[u_1, u_2]\times [v_1, v_2]$.

\noindent Type 1: $R$ is above the curve of $\phi$, \emph{i.e.}, $R\subset\{(u, v)\mid v\geq\phi(u)\}$.

For this type of rectangle, there exist $s_1, s_2\in [0, 1]$ such that $v_1=\phi(s_1)$ and $v_2=\phi(s_2)$. Clearly, $u_1<u_2\leq s_1<s_2$.
We distinguish three cases to show that $V(R)\geq 0$.
\begin{itemize}%\setlength{\itemindent}{5pt}
\item[\rm (a)] $A_{\Gamma_\phi}(u_1, \phi(s_1))=u_1\leq \phi(s_1)-\max\limits_{t\in [u_1, s_1]}\{\phi(t)-\Gamma_\phi(t)\}$.

In this case,  it follows from the increasingness of $A_{\Gamma_\phi}$ that $A_{\Gamma_\phi}(u_1, \phi(s_2))=u_1$, and so
$V(R)=A_{\Gamma_\phi}(u_2, \phi(s_2))-A_{\Gamma_\phi}(u_2, \phi(s_1))\geq 0$.

\item[\rm (b)] $A_{\Gamma_\phi}(u_1, \phi(s_1))= \phi(s_1)-\max\limits_{t\in [u_1, s_1]}\{\phi(t)-\Gamma_\phi(t)\}<u_1$ and

$A_{\Gamma_\phi}(u_2, \phi(s_2))=u_2\leq \phi(s_2)-\max\limits_{t\in [u_2, s_2]}\{\phi(t)-\Gamma_\phi(t)\}$.

In this case, we have
\begin{eqnarray*}
V(R)&=& \phi(s_1)-\max\limits_{t\in [u_1, s_1]}\{\phi(t)-\Gamma_\phi(t)\}-A_{\Gamma_\phi}(u_2, \phi(s_1))-A_{\Gamma_\phi}(u_1, \phi(s_2))+u_2\\
&\geq& \max\limits_{t\in [u_2, s_1]}\{\phi(t)-\Gamma_\phi(t)\}-\max\limits_{t\in [u_1, s_1]}\{\phi(t)-\Gamma_\phi(t)\}-u_1+u_2\\
&=&\min\left\{\max\limits_{t\in [u_2, s_1]}\{\phi(t)-\Gamma_\phi(t)\}-\max\limits_{t\in [u_1, u_2]}\{\phi(t)-\Gamma_\phi(t)\}, 0\right\}-u_1+u_2\\
&\geq&\min\left\{u_2-u_1+\phi(u_2)-\Gamma_\phi(u_2)-\max\limits_{t\in [u_1, u_2]}\{\phi(t)-\Gamma_\phi(t)\}, u_2-u_1\right\}.
\end{eqnarray*}
For all $t\in [u_1, u_2]$, it holds that
\begin{eqnarray*}
 & &  u_2-u_1+\phi(u_2)-\Gamma_\phi(u_2)-(\phi(t)-\Gamma_\phi(t))\\
 &=& (t-u_1)+(\Gamma_\phi(t)-t-\phi(t))-(\Gamma_\phi(u_2)-u_2-\phi(u_2)) \\
 &\geq& 0.
\end{eqnarray*}
Hence, we have $V(R)\geq 0$.

\item[\rm (c)] $A_{\Gamma_\phi}(u_1, \phi(s_1))= \phi(s_1)-\max\limits_{t\in [u_1, s_1]}\{\phi(t)-\Gamma_\phi(t)\}<u_1$ and

$A_{\Gamma_\phi}(u_2, \phi(s_2))=\phi(s_2)-\max\limits_{t\in [u_2, s_2]}\{\phi(t)-\Gamma_\phi(t)\}<u_2$.

In this case, it follows from the increasingness of $A_{\Gamma_\phi}$ that $A_{\Gamma_\phi}(u_2, \phi(s_1))=\phi(s_1)-\max\limits_{t\in [u_2, s_1]}\{\phi(t)-\Gamma_\phi(t)\}$, and so
\begin{eqnarray*}
V(R) &=& \max\limits_{t\in [u_2, s_1]}\{\phi(t)-\Gamma_\phi(t)\}-\max\limits_{t\in [u_1, s_1]}\{\phi(t)-\Gamma_\phi(t)\}
-\max\limits_{t\in [u_2, s_2]}\{\phi(t)-\Gamma_\phi(t)\}+\phi(s_2)\\
& &-\min\{u_1,  \phi(s_2)-\max\limits_{t\in [u_1, s_2]}\{\phi(t)-\Gamma_\phi(t)\}\}.
\end{eqnarray*}
We further distinguish two subcases.
\end{itemize}

\begin{itemize}%\setlength{\itemindent}{5pt}
\item[\rm (c1)] $\max\limits_{t\in [u_2, s_1]}\{\phi(t)-\Gamma_\phi(t)\}\geq \min\{\max\limits_{t\in [u_1, u_2]}\{\phi(t)-\Gamma_\phi(t)\},
\max\limits_{t\in [s_1, s_2]}\{\phi(t)-\Gamma_\phi(t)\}\}$.

In this subcase, we have
\[\max\limits_{t\in [u_2, s_1]}\{\phi(t)-\Gamma_\phi(t)\}= \min\{\max\limits_{t\in [u_1, s_1]}\{\phi(t)-\Gamma_\phi(t)\},
\max\limits_{t\in [u_2, s_2]}\{\phi(t)-\Gamma_\phi(t)\}\}.\]
Note that
\[\max\limits_{t\in [u_1, s_2]}\{\phi(t)-\Gamma_\phi(t)\}= \max\{\max\limits_{t\in [u_1, s_1]}\{\phi(t)-\Gamma_\phi(t)\},
\max\limits_{t\in [u_2, s_2]}\{\phi(t)-\Gamma_\phi(t)\}\}.\]
Hence,
\begin{eqnarray*}
 & & \max\limits_{t\in [u_2, s_1]}\{\phi(t)-\Gamma_\phi(t)\}-\max\limits_{t\in [u_1, s_1]}\{\phi(t)-\Gamma_\phi(t)\}
-\max\limits_{t\in [u_2, s_2]}\{\phi(t)-\Gamma_\phi(t)\}+\phi(s_2)\\
& &-(\phi(s_2)-\max\limits_{t\in [u_1, s_2]}\{\phi(t)-\Gamma_\phi(t)\})\\
&=& \max\limits_{t\in [u_2, s_1]}\{\phi(t)-\Gamma_\phi(t)\}-\max\limits_{t\in [u_1, s_1]}\{\phi(t)-\Gamma_\phi(t)\}
-\max\limits_{t\in [u_2, s_2]}\{\phi(t)-\Gamma_\phi(t)\}\\
& &+\max\limits_{t\in [u_1, s_2]}\{\phi(t)-\Gamma_\phi(t)\}\\
&=& 0,
\end{eqnarray*}
and so $V(R)\geq 0$.

\item[\rm (c2)] $\max\limits_{t\in [u_2, s_1]}\{\phi(t)-\Gamma_\phi(t)\}< \min\{\max\limits_{t\in [u_1, u_2]}\{\phi(t)-\Gamma_\phi(t)\},
\max\limits_{t\in [s_1, s_2]}\{\phi(t)-\Gamma_\phi(t)\}\}$.

In this subcase, there exist $u_0\in [u_1, u_2[\,$ and $s_0\in \,]s_1, s_2]$ such that
\[\max\limits_{t\in [u_1, s_1]}\{\phi(t)-\Gamma_\phi(t)\}=\phi(u_0)-\Gamma_\phi(u_0)>0\]
and
\[\max\limits_{t\in [u_2, s_2]}\{\phi(t)-\Gamma_\phi(t)\}=\phi(s_0)-\Gamma_\phi(s_0)>0,\]
and so
\[\phi(u_0)-\Gamma_\phi(u_0)>\max\{0, \phi(s_1)-u_1\}\]
and \[\phi(s_0)-\Gamma_\phi(s_0)>\max\{0, \phi(s_2)-u_2\}.\]
It is easy to see that $\Gamma_\phi(u_0)<m_\phi(u_0)$ and $\Gamma_\phi(s_0)<m_\phi(s_0)$, and thus there exist $i_0, j_0\in I$
such that $u_0\in\,]a_{i_0}, b_{i_0}[\,$ and $s_0\in\,]a_{j_0}, b_{j_0}[\,$. Since $\phi(t)-\Gamma_\phi(t)$ is strictly increasing on $[a_i, u_i^\ast]$ and is strictly decreasing on $[u_i^\ast, b_i]$ for all $i\in I$, we have $i_0\not=j_0$, and thus $a_{j_0}\geq b_{i_0}$.
Note that
\begin{eqnarray*}
V(R)&\geq & \max\limits_{t\in [u_2, s_1]}\{\phi(t)-\Gamma_\phi(t)\}-\max\limits_{t\in [u_1, s_1]}\{\phi(t)-\Gamma_\phi(t)\}
-\max\limits_{t\in [u_2, s_2]}\{\phi(t)-\Gamma_\phi(t)\}+\phi(s_2)-u_1\\
&\geq& -(\phi(u_0)-\Gamma_\phi(u_0))-(\phi(s_0)-\Gamma_\phi(s_0))+\phi(s_2)-u_1\\
&\geq& -(\phi(u_0)-\Gamma_\phi(u_0))-(\phi(s_0)-\Gamma_\phi(s_0))+\phi(s_0)-u_0\\
&=& \Gamma_\phi(s_0)+\Gamma_\phi(u_0)-u_0-\phi(u_0)\\
&\geq & \Gamma_\phi(a_{j_0})+\Gamma_\phi(b_{i_0})-b_{i_0}-\phi(b_{i_0})\\
&= & m_\phi(a_{j_0})-\max\{b_{i_0}, \phi(b_{i_0})\}.
\end{eqnarray*}
Hence, it follows from condition (ii) that $V(R)\geq 0$.
\end{itemize}

\noindent Type 2: $R$ is below the curve of $\phi$, \emph{i.e.}, $R\subset\{(u, v)\mid v\leq\phi(u)\}$.

The proof is similar to that of Type 1.

\noindent Type 3: $R$ has two vertices on the curve of $\phi$, \emph{i.e.}, $R=[u_1, u_2]\times [\phi(u_1), \phi(u_2)]$.

For this type of rectangle, we have
\begin{eqnarray*}
V(R) &=& \Gamma_\phi(u_1)+ \Gamma_\phi(u_2)-\min\{u_1, \phi(u_2)-\max_{t\in [u_1, u_2]}\{\phi(t)-\Gamma_\phi(t)\}\}\\
& &-\min\{\phi(u_1), u_2-\max_{t\in [u_1, u_2]}\{t-\Gamma_\phi(t)\}\}\\
&=&\max\{\Gamma_\phi(u_2)+ \Gamma_\phi(u_1)-u_1-\phi(u_1), \Gamma_\phi(u_1)+ \Gamma_\phi(u_2)-u_1-u_2+\max_{t\in [u_1, u_2]}\{t-\Gamma_\phi(t)\},\\
& &\Gamma_\phi(u_1)+ \Gamma_\phi(u_2)-\phi(u_1)-\phi(u_2)+\max_{t\in [u_1, u_2]}\{\phi(t)-\Gamma_\phi(t)\},\\
& &\Gamma_\phi(u_1)+ \Gamma_\phi(u_2)-u_2-\phi(u_2)+\max_{t\in [u_1, u_2]}\{t-\Gamma_\phi(t)\}+\max_{t\in [u_1, u_2]}\{\phi(t)-\Gamma_\phi(t)\}\}.
\end{eqnarray*}
We distinguish three cases to show that $V(R)\geq 0$.
\begin{itemize}%\setlength{\itemindent}{5pt}
\item[\rm (a)] $\Gamma_\phi(u_1)=m_\phi(u_1)$ or $\Gamma_\phi(u_2)=m_\phi(u_2)$.

In this case, since
\begin{eqnarray*}
 & &  \Gamma_\phi(u_1)+ \Gamma_\phi(u_2)-u_1-u_2+\max_{t\in [u_1, u_2]}\{t-\Gamma_\phi(t)\}\\
 &\geq& \max\{\Gamma_\phi(u_1)-u_1, \Gamma_\phi(u_2)-u_2\}
\end{eqnarray*}
and
\begin{eqnarray*}
 & &  \Gamma_\phi(u_1)+ \Gamma_\phi(u_2)-\phi(u_1)-\phi(u_2)+\max_{t\in [u_1, u_2]}\{\phi(t)-\Gamma_\phi(t)\}\\
 &\geq& \max\{\Gamma_\phi(u_1)-\phi(u_1), \Gamma_\phi(u_2)-\phi(u_2)\},
\end{eqnarray*}
we have
\begin{eqnarray*}
V(R) &\geq& \max\{\Gamma_\phi(u_1)-u_1, \Gamma_\phi(u_1)-\phi(u_1), \Gamma_\phi(u_2)-u_2, \Gamma_\phi(u_2)-\phi(u_2)\}\\
&=& \max\{\Gamma_\phi(u_1)-m_\phi(u_1), \Gamma_\phi(u_2)-m_\phi(u_2)\}\\
&=& 0.
\end{eqnarray*}

\item[\rm (b)] There exists $i\in I$ such that $\{u_1, u_2\}\subset\,]a_i, b_i[\,$.

In this case, since $t-\Gamma_\phi(t)$ and $\phi(t)-\Gamma_\phi(t)$ are strictly increasing on $[a_i, u_i^\ast]$ and are strictly decreasing on $[u_i^\ast, b_i]$, if $u_i^\ast\in[u_1, u_2]$, then
\begin{eqnarray*}
 & &  \Gamma_\phi(u_1)+ \Gamma_\phi(u_2)-u_2-\phi(u_2)+\max_{t\in [u_1, u_2]}\{t-\Gamma_\phi(t)\}+\max_{t\in [u_1, u_2]}\{\phi(t)-\Gamma_\phi(t)\}\\
 &=& \Gamma_\phi(u_1)+ \Gamma_\phi(u_2)-u_2-\phi(u_2)+u_i^\ast-\Gamma_\phi(u_i^\ast)+\phi(u_i^\ast)-\Gamma_\phi(u_i^\ast) \\
 &=& (\Gamma_\phi(u_1)-\Gamma_\phi(u_i^\ast))+ (\Gamma_\phi(u_2)-u_2-\phi(u_2))-(\Gamma_\phi(u_i^\ast)-u_i^\ast-\phi(u_i^\ast)) \\
 &=& 0.
\end{eqnarray*}
If $u_i^\ast\notin[u_1, u_2]$, then
\begin{eqnarray*}
 & &  \Gamma_\phi(u_1)+ \Gamma_\phi(u_2)-u_2-\phi(u_2)+\max_{t\in [u_1, u_2]}\{t-\Gamma_\phi(t)\}+\max_{t\in [u_1, u_2]}\{\phi(t)-\Gamma_\phi(t)\}\\
 &=& \max\{\Gamma_\phi(u_1)-\Gamma_\phi(u_2), (\Gamma_\phi(u_2)-u_2-\phi(u_2))-(\Gamma_\phi(u_1)-u_1-\phi(u_1))\} \\
 &=& 0.
\end{eqnarray*}
Hence, we have $V(R)\geq 0$ in both subcases.

\item[\rm (c)] There exist $i, j\in I$ with $i\not=j$ such that $u_1\in\,]a_i, b_i[\,$ and $u_2\in\,]a_j, b_j[\,$.

In this case, we have $a_j\geq b_i$ and
\begin{eqnarray*}
& & \Gamma_\phi(u_2)+ \Gamma_\phi(u_1)-u_1-\phi(u_1) \\
&\geq&  \Gamma_\phi(a_j)+ \Gamma_\phi(b_i)-b_i-\phi(b_i) \\
&=& m_\phi(a_j)-\max\{b_i, \phi(b_i)\}.
\end{eqnarray*}
Hence, it follows from condition (ii) that $V(R)\geq 0$.
\end{itemize}
We conclude that $A_{\Gamma_\phi}$ is a copula.
\end{proof}

\begin{remk}
For any $\phi\in\Phi$, we have the following observations:
\begin{itemize}  %\setlength{\itemindent}{5pt}
\item[{\rm(i)}] if $\Gamma_\phi\in\Delta_\phi$ is determined by pairwisely disjoint countably many open subintervals $\{\,]a_i, b_i[\,\}_{i\in I}$ of $\,]0, 1[\,$,
then it is not strictly increasing.
Since the curvilinear section $\Gamma_\phi$ of $\Pi$ given by $\Gamma_\phi(t)=t\phi(t)$ is strictly increasing, $A_{\Gamma_\phi}$ is not a copula;

\item[{\rm(ii)}] $\phi$ is trivially compatible with any single open subinterval of $\,]0, 1[\,$.
Since the curvilinear section $\Gamma_\phi$ of $W$ given by $\Gamma_\phi(t)=\max\{0, t+\phi(t)-1\}$
is determined by the single open interval $\{\,]0, 1[\,\}$, $A_{\Gamma_\phi}$ is a copula given by
\[
A_{\Gamma_\phi}(u, v) =\left \{
        \begin {array}{ll}
         \min\{v, u-t^\ast\}
                  &\quad \text{if\ \ $\phi^{-1}(v)<t^\ast<u$}  \\[2mm]
         \min\{u, v-\phi(t^\ast)\}
                  &\quad \text{if\ \ $u<t^\ast<\phi^{-1}(v)$}\\[2mm]
         W(u, v)
                  &\quad \text{otherwise,}
        \end {array}
       \right.
\]
where $t^\ast\in\,]0, 1[\,$ is the unique solution to the equation: $t+\phi(t)-1=0$.
Interestingly, if $t^\ast=\dfrac{1}{2}$, \emph{i.e.}, $\dfrac{1}{2}$ is a fixed point of $\phi$, then $A_{\Gamma_\phi}=A_\delta$,
where $\delta$ is the diagonal section of $W$ given by $\delta(t)=\max\{0, 2t-1\}$;

\item[{\rm(iii)}] there exist pairwisely disjoint countably infinite open subintervals of $\,]0, 1[\,$ compatible with $\phi$. Explicitly, starting from any open interval
$\,]a_1, b_1[\,\subseteq \,]0, 1[\,$ with $b_1<1$, for any positive integer $n$, define $a_{n+1}$ and $b_{n+1}$ as follows:
\[a_{n+1}=\max\{\phi(b_n), \phi^{-1}(b_n)\};\ \ b_{n+1}=\frac{1+a_{n+1}}{2}.\]
It is easy to see that $\phi$ is compatible with $\{\,]a_n, b_n[\,\}_{n=1}^\infty$.
\end{itemize}
\end{remk}

\section{Concluding remarks}\label{Sec-Conclu}

We have characterized when $A_{\Gamma_\phi}$ is a copula (see Theorem~\ref{Thrm-main}). For future work, it is interesting to consider the following problem: \\
For any $\phi\in\Phi$ and $\Gamma_\phi\in\Delta_\phi$, denote by $\mathcal C_{\Gamma_\phi}$ the set of all copulas
with curvilinear section $\Gamma_\phi$. Note that $B_{\Gamma_\phi}$ defined by~\eqref{Eq-DfntofBGamma} belongs to
$\mathcal C_{\Gamma_\phi}$. Let $\overline{C}_{\Gamma_\phi}$ be the supremum of $\mathcal C_{\Gamma_\phi}$, \emph{i.e.},
$\overline{C}_{\Gamma_\phi}(u, v)=\sup\{C(u, v)\mid C\in\mathcal C_{\Gamma_\phi}\}$. Clearly, $\overline{C}_{\Gamma_\phi}\leq A_{\Gamma_\phi}$.
It is known that $\overline{C}_{\Gamma_\phi}$ is a quasi-copula (see~\cite{NQRU04}) with curvilinear section $\Gamma_\phi$, but it may not be a copula
even if $\phi={\rm id}$ (see~\cite{NQRU08}). However, if $A_{\Gamma_\phi}$ is a copula, then $\overline{C}_{\Gamma_\phi}$ is also a copula as
$\overline{C}_{\Gamma_\phi}=A_{\Gamma_\phi}$. Naturally, the problem is to characterize when $\overline{C}_{\Gamma_\phi}$ is a copula.

\end{document}